\documentclass[10pt]{amsart}

\def\cal{\mathcal}

\let\newpf\proof \let\proof\relax

\def\0{{\mathbf{0}}}

\def\g{\gamma}

\def\be{\begin{equation}}
\def\ee{\end{equation}}

\def\d{{\underline d}}

\newtheorem{thm}{Theorem}[section]

\newtheorem{lemma}[thm]{Lemma}

\theoremstyle{remark}
\newtheorem{rem}{Remark}[section]

\numberwithin{equation}{section}

\def \bn {\hfill \\ \smallskip\noindent}

\theoremstyle{definition}

\def\proof{\bn {\bf Proof.} }

\newcommand{\inter}{\operatorname{int}}

\newcommand{\CC}{{\cal C}}
\newcommand{\DD}{{\cal D}}

\newcommand{\II}{{\cal I}}
\newcommand{\FF}{{\cal F}}

\newcommand{\QQ}{{\cal Q}}

\newcommand{\C}{{\mathbb C}}

\newcommand{\N}{{\mathbb N}}

\newcommand{\R}{{\mathbb R}}

\newcommand{\Z}{{\mathbb Z}}

\def\B0{{\bold{0}}}

\catcode`\@=12

\def\Empty{}
\newcommand\oplabel[1]{
  \def\OpArg{#1} \ifx \OpArg\Empty {} \else
        \label{#1}
  \fi}
                
\newcommand{\comm}[1]{}
\newcommand{\comment}[1]{}

\begin{document}

\title{Bifurcations of unimodal maps}

\author{Artur Avila and Carlos Gustavo Moreira}

\address{
Coll\`ege de France -- 3 Rue d'Ulm \\
75005 Paris -- France.
}
\email{avila@impa.br}

\address{
IMPA -- Estr. D. Castorina 110 \\
22460-320 Rio de Janeiro -- Brazil.
}
\email{gugu@impa.br}

\thanks{Partially supported by Faperj and CNPq, Brazil.}

\date{\today}

\begin{abstract}

We review recent results that lead to a very precise understanding of the
dynamics of typical unimodal maps from the statistical point of view. 
We also describe the (generalized) renormalization approach
to the study of the statistical properties of typical unimodal maps.

\end{abstract}

\setcounter{tocdepth}{1}

\maketitle


\section{Introduction}

A unimodal map is a smooth (at least $C^2$) map $f:I \to I$  of an interval
$I \subset \R$ with a unique critical point $c \in \inter I$ which is a
maximum.  We will also assume that $f(\partial I) \subset \partial I$.
The main examples of unimodal maps are given by the quadratic
family of maps $p_a(x)=a-x^2$, where $-1/4 \leq a \leq 2$ is a parameter.

In this work we will describe in detail the properties of
``typical'' unimodal maps.  Here typical is to be understood in the
measure-theoretical sense: it should correspond to a generalization
of ``almost every parameter'' for the quadratic family.

\subsection{Regular maps}

A unimodal map is said to be {\it Kupka-Smale} if its critical point is
non-degenerate and if it has only hyperbolic periodic orbits.

A unimodal map is said to be {\it hyperbolic} if there are finitely many
hyperbolic periodic sinks and the dynamical
interval can be written as a union $I=U \cup K$ into invariant sets where
$f|K$ is uniformly expanding and every $x \in U$ is attracted by some
hyperbolic periodic sink.  In this case, $K$ has always Lebesgue measure
zero.  Thus hyperbolic maps have deterministic dynamics.

A unimodal map is said to be {\it regular} if it is Kupka-Smale, and
if its critical point is attracted to a periodic sink and is not
periodic or preperiodic.

By a result of Ma\~n\'e, a regular unimodal map is always hyperbolic, and
it follows that the set of regular maps is open.  The
work of Lyubich, Graczyk-Swiatek and Kozlovski shows that regular maps are
dense among quadratic/smooth/analytic unimodal maps.
It is easy to see that regular maps are {\it structurally stable}: any
$C^2$ small perturbation is still a unimodal map and is
topologically conjugate to the original map.  Moreover, the converse also
holds: if $f$ is structurally stable among quadratic/smooth/analytic
maps then $f$ is regular.

Let us say that an analytic family of unimodal maps is {\it non-degenerate}
if regular parameters are dense.  The quadratic family is an example of a
non-degenerate family, and it is possible to show that non-degeneracy is a
very weak assumption: for instance, any analytic family of
unimodal maps with {\it negative Schwarzian derivative}\footnote
{That is, $Sf=D^3 f-\frac {3} {2}
(D^2 f)^2<0$.  This condition defines an open set of maps which contains the
quadratic family.} is non-degenerate provided it contains one
regular parameter (in particular, any analytic family $C^3$ close to the
quadratic family is non-degenerate).  Throughout this paper, a property will
be said to be typical for analytic unimodal maps if it is satisfied for
almost every parameter in any non-degenerate analytic family of
unimodal maps.  For smooth ($C^k$, $k=2,3,...,\infty$) unimodal maps,
typical will denote ``for almost every parameter in any family belonging
to some generic (residual) set of $C^k$ families of unimodal maps''.
For instance, Kupka-Smale maps are typical.

Although regular maps are topologically generic, their complement is
certainly non-negligible: in the quadratic family, and actually in any
family $C^2$ close to the quadratic family, the set of non-regular
parameters has positive Lebesgue measure (a version of
Jakobson's Theorem, see for instance \cite {T}).
Thus, in our description of the
dynamics of typical unimodal maps, non-regular maps must be present (and
actually form the interesting part of the description).

\subsection{Typical non-regular maps: topological description}

Let us say that $T \subset I$ is a restrictive interval (of period $m$) if
it contains the critical point $c$, $f^m(T) \subset T$, $f^m(\partial T)
\subset \partial T$, and $f^j(T) \cap
\inter T=\emptyset$, $1 \leq j \leq m-1$ (if $m>1$ then $T$ is also called a
{\it renormalization interval}).  Notice that $f^m:T \to T$
is unimodal.  Letting $\hat T=[f^{2m}(c),f^m(c)]$, it is easy to see that
either $c \in \inter \hat T$ and $f(\hat T)=\hat T$ or $f^m|T$ has trivial
dynamics: all orbits are asymptotic to a periodic orbit.  The following
result was proved in this generality in \cite {AM4}, building on the key
quadratic case which was proved by Lyubich \cite {regular} and on the work
of \cite {ALM}.

\begin{thm} \label {top}

Let $f$ be a typical non-regular unimodal map.  Then:

\begin{enumerate}

\item $f$ is Kupka-Smale and has finitely many hyperbolic periodic sinks,

\item $f$ has a smallest restrictive interval $T$ and the first return map
$f^m|T$ is conjugate to a quadratic map,

\item There is a decomposition $I=U \cup K$ in invariant sets where
$U$ is the set of points $x \in I$ which are either attracted to sinks or are
eventually trapped in $\hat T$, and $f|K$ is uniformly expanding.

\item $f^m:\hat T \to \hat T$ is topologically mixing.

\end{enumerate}

\end{thm}

Thus one sees that to understand a typical non-regular unimodal map,
one has to understand the dynamics of orbits in the ``attractor''
$A=\cup_{k=0}^{m-1} f^k(\hat T)$, as the complementary dynamics is
hyperbolic.

It turns out that, in a certain sense, the attractor $A$ can not be
decomposed further: it is a genuine topological and metric attractor in the
sense of Milnor: the set of points $x \in A$
whose orbit is dense in $A$ is both residual and a full measure
subset of $A$
(this actually holds under much milder
assumption than ``typical'', see \cite {attractors}).

\subsection{Statistical description}

In order to describe the dynamics of a typical non-regular map from the
statistical point of view, we consider invariant probability
measures $\mu$.  We say
that $x \in I$ is in the basin of $\mu$ if for any continuous $\phi:I \to
\R$ we have:
\be
\lim \frac {1} {m} \sum_{k=0}^{m-1} \phi(f^k(x))=\int \phi d\mu.
\ee
We are interested in measures $\mu$ which possess a basin of positive
Lebesgue measure.  Such measures will be called physical.  Examples of
physical measures are given by the invariant measures supported on
hyperbolic periodic sinks.  Other examples are provided by ergodic
absolutely continuous invariant probabilities (they are physical measures by
Birkhoff's Ergodic Theorem).

For a typical non-regular map one clearly has (using Theorem \ref {top})
a finite (possibly zero)
number of physical measures corresponding to sinks and possibly
other physical measures
supported on $A$.  It turns out that $f|A$ is ergodic with respect to
Lebesgue measure, so there can be at most one physical measure supported in
$A$, and if it exists then its basin must have
full Lebesgue measure in $A$.
A particularly nice situation occurs if $f$ has an
ergodic invariant measure equivalent to Lebesgue measure on
$A$ (in this case we will say that $f$ is stochastic).

It turns out that there are quadratic maps (with an attractor $A$ as above)
without physical measures \cite {Jo}.  There are also examples of
quadratic maps with rather unexpected physical measures: in \cite {HK}
it is shown that the physical measure can be supported on a
hyperbolic repelling fixed point (which might be called the
``statistical attractor'' of the map).  So one is naturally led to
ask if those possibilities really arise in the typical setting.  The
following result was proved in \cite {ALM} and \cite {AM4}, and is based on
the work of Lyubich \cite {regular} which covers the quadratic case.

\begin{thm}

A typical non-regular unimodal map possesses a unique ergodic
invariant measure $\mu$ equivalent to Lebesgue measure on $A$.

\end{thm}
 
To describe further properties of $\mu$, it is convenient to consider the
system $f^m|\hat T$ which has an absolutely continuous invariant measure
$\hat \mu=m \mu|\hat T$.  The system $(f^m|\hat T,\hat \mu)$ is always
mixing (actually, a general result of Ledrappier
implies that $(f^m|\hat T,\hat \mu)$ is weak Bernoulli) and has a positive
Lyapunov exponent.  The following result was proved in \cite {AM2}
and \cite {AM4} (the quadratic case was covered in \cite {AM1}).

\begin{thm} \label {stat}

For a typical non-regular unimodal map, $(f^m|\hat T,\hat \mu)$ has
exponential decay of correlations and is stochastically stable.

\end{thm}

The proof of the previous result is based on Theorem \ref {top} coupled with
a good description of the critical orbit in terms of hyperbolicity and
recurrence.  Indeed, Keller-Nowicki \cite {KN} and Young \cite {Y}
have shown that if $f$ satisfies the {\it Collet-Eckmann condition},
that is, $|Df^n(f(c))|>C \lambda^n$, for some $C>0$, $\lambda>1$,
then the system $(f^m|\hat T,\hat
\mu)$ has exponential decay of correlations
provided that $f$ is Kupka-Smale and $f^m|\hat T$ is topologically
mixing\footnote {It turns out that under those condition, Collet-Eckmann is
actually equivalent to exponential decay of correlations, see \cite {NS}.}
(due to Theorem \ref {top}, those conditions hold for typical maps).  It has
been shown by Baladi-Viana \cite {BV} that, under the additional hypothesis
of {\it subexponential recurrence} of the critical orbit,
that is $|f^n(c)-c|>e^{-\alpha n}$
for every $\alpha>0$ and for every $n$ sufficiently big, the system
$(f^m|\hat T,\hat \mu)$ is stochastically stable (they actually
assume that $f$ is at least $C^3$, for the $C^2$ case one must use a
result of Tsujii \cite {T2}).  Thus, Theorem \ref {stat} is actually a
consequence of the following:

\begin{thm}

A typical non-regular unimodal map is Collet-Eckmann and the recurrence of
its critical orbit is subexponential.

\end{thm}

Better estimates can be given for typical analytic unimodal maps: in this
case the recurrence of the critical orbit is actually polynomial with
exponent $1$, that is,
\be
\limsup \frac {-\ln |f^n(c)-c|} {\ln n}=1.
\ee

\subsection{Unimodal maps from the point of view of
(generalized) renormalization}

It is clear that the complications in the study of the dynamics of
unimodal maps arise from the presence of the critical point.  In a gross
simplification, one can identify two approaches to face the complications
posed by the critical point:

\begin{enumerate}

\item Focus on the good part, that is, concentrate
on the description of the dynamics in large scales and
away from the critical point,

\item Focus on the problematic part, that is, concentrate on the description
of the dynamics in small scales and near the critical point.

\end{enumerate}

The best example of the first approach is the {\it inducing method}.  This
method was used by Jakobson to obtain a positive measure set of parameters
in the quadratic family for which he constructed absolutely continuous
invariant measures.  In his work, he constructs an induced Markov map
by finding more and more branches of iterates of $f$ which reach large
scale: in the end he obtains a partition (modulo 0) of the phase space in
intervals $T^j$ such that a suitable iterate
$f^{r_j}|T^j$ is a diffeomorphism over an interval
of definite size.  To give an example of a more recent application,
convenient induced Markov maps can also be used to obtain fine
statistical properties of unimodal (and multimodal) maps under convenient
assumptions on the critical behavior \cite {BLS}.

The best example of the second approach is the {\it renormalization method}. 
According to the description of \cite {cambridge}, in this method
one considers a sequence of small intervals around the critical
point and looks at their first return maps, which are also called
generalized renormalizations of the initial system.  One application of this
method is to show that, under certain combinatorial assumptions,
the geometry of the critical orbit is rigid.

Of course, both approaches are not completely separated, but their
philosophy is quite distinct.  Here we will adopt the renormalization point
of view and use it to describe the dynamics of unimodal maps in all
respects.  We will start with some combinatorial preparation, then describe
the Phase-Parameter relation, and finally we will discuss the statistical
arguments involved.  We will focus on the quadratic case for simplicity, and
our presentation can be seen as an informal guide to \cite {AM1}.

\comm{
and then discuss
the non-trivial parts of the generalization of those ideas for the
analytic case, mainly related to the Phase-Parameter relation
(the smooth case follows from the analytic case by rather
different considerations).
}

\section{Statistical properties of the quadratic family}

Let us normalize the quadratic family as
\be
f_a(x)=a-1-a x^2,
\ee
where $1/2 \leq a \leq 2$, so that $p_a$ is a unimodal map in the canonical
interval $I=[-1,1]$.

\subsection{Combinatorics and the phase-parameter relation}

\subsubsection{Renormalization}

We say that $f$ is {\it renormalizable} if there is an interval $0 \in T$
and $m>1$ such that $f^m(T) \subset T$ and
$f^j(\inter T) \cap \inter T=\emptyset$ for $1 \leq j<m$.
The maximal such interval is called the {\it renormalization interval of
period $m$}, it has the property that $f^m(\partial T) \subset \partial T$.

The set of renormalization periods of $f$ gives an increasing (possibly
empty) sequence of numbers $m_i$, $i=1,2,...$,
each related to a unique renormalization
interval $T^{(i)}$ which form a nested sequence of intervals.  We include
$m_0=1$, $T^{(0)}=I$ in the sequence to simplify the notation.

We say that $f$ is {\it finitely renormalizable} if there is a smallest
renormalization interval $T^{(k)}$.  We say that
$f \in \FF$ if $f$ is finitely
renormalizable and $0$ is recurrent but not periodic.  We let $\FF_k$ denote
the set of maps $f$ in $\FF$ which are exactly $k$ times renormalizable.

The analysis of infinitely renormalizable maps is quite different from the
finitely renormalizable case.  The following fundamental
result was proved in \cite {regular}:

\begin{thm}

The set of infinitely renormalizable parameters has zero Lebesgue measure in
the quadratic family.

\end{thm}

The proof of this result goes beyond the scope of this note (see
the survey \cite {notices} for a discussion of some of the elements of the
proof).

A much simpler argument shows that almost every quadratic map with a
non-recurrent critical point is indeed regular.  Thus, we may concentrate on
quadratic maps $f \in \FF$.

\subsubsection{Principal nest}

We say that a symmetric interval $T \subset I$ is nice for $f$ if
$f^n(\partial T) \cap \inter T=\emptyset$, $n \geq 0$.  It is easy to see
that the first return map $R_T$ to a nice interval $T$ has the following
property: its domain is a disjoint union of intervals $T^j$ and $R_T|T^j$
is a diffeomorphism onto $T$ if $0 \notin T^j$.  On the other hand, the
component of $0$ (if it exists) is itself a nice interval.  We will now
introduce a special sequence of nice intervals obtained by iteration of this
procedure, which will be the basis of our analysis of finitely
renormalizable maps (this sequence is also very important in the analysis of
infinitely renormalizable maps as well).

Let $\Delta_k$ denote the set of all maps $f$
which have (at least) $k$ renormalizations and which have an orientation
reversing non-attracting periodic point of period $m_k$ which we denote
$p_k$ (that is, $p_k$ is the fixed point of
$f^{m_k}|_{T^{(k)}}$ with $Df^{m_k}(p_k) \leq -1$).
For $f \in \Delta_k$, we denote $T^{(k)}_0=[-p_k,p_k]$.
We define by induction a (possibly
finite) sequence $T^{(k)}_i$, such that $T^{(k)}_{i+1}$ is the component of
the domain of $R_{T^{(k)}_i}$ containing $0$.  If this sequence is
infinite, then either it converges to a point or to an interval.

If $\cap_i T^{(k)}_i$ is a point, then
$f$ has a recurrent critical point which is not
periodic, and it is possible to show that $f$ is not $k+1$ times
renormalizable.  Obviously in this case we have
$f \in \FF_k$, and all maps in $\FF_k$ are obtained in this way:
if $\cap_i T^{(k)}_i$ is an
interval, it is possible to show that $f$ is $k+1$ times renormalizable.

It is important to notice that the domain of the first return map
to $T^{(k)}_i$ is always dense in $T^{(k)}_i$.  Moreover, the next result
shows that, outside a very special case,
the return map has a hyperbolic structure.

\begin{lemma} \label {hyperbol}

Assume that $T^{(k)}_i$ does not have a non-hyperbolic periodic orbit
in its boundary.  For all $T^{(k)}_i$ there exists $C>0$, $\lambda>1$
such that if $x,f(x),...,f^{n-1}(x)$ do not belong to
$T^{(k)}_i$ then $|Df^n(x)|>C \lambda^n$.

\end{lemma}

Since almost every non-regular map belongs to $\FF_\kappa$ for some $\kappa$,
it is enough to work with non-regular maps in some fixed $\FF_\kappa$.  Once
$\kappa$ is fixed, we may introduce the following convenient notation for $f
\in \Delta_\kappa$.  

Let $I_n=T^{(\kappa)}_n$, and let $R_n$ be the first return map to $I_n$. 
The domain of $R_n$ is a union of intervals (labeled by a subset of $\Z$)
denoted $I^j_n$, where we reserve the index $0$ for the component of the
critical point: $0 \in I^0_n=I_{n+1}$.  Since $I_n$ is nice,
$R_n|I^j_n$ is a diffeomorphism onto $I_n$ whenever $j \neq 0$.  The return
to level $n$ will be called central if $R_n(0) \in I_{n+1}$.

Let $\Omega$ be the set of all finite words $\d=(j_1,...,j_m)$ of non-zero
integers, and let $|\d|=m$ denote the length of $\d$.
If $\d=(j_1,...,j_m) \in \Omega$,
let $I^\d_n=\{x \in I_n,\, R_n^{k-1}(x) \in I^{j_i}_n,\, 1 \leq k \leq m\}$,
and define $R^\d_n:I^\d_n \to I_n$ by $R^\d_n=R^m_n$.  Let
$C^\d_n=(R^\d_n)^{-1}(I_{n+1})$.  The map $L_n:\cup_{\d \in \Omega} C^\d_n
\to I_{n+1}$, $L_n|C^\d_n=R^\d_n$ is the first landing map from $I_n$ to
$I_{n+1}$.

The dependence on $f$ is implicit in the above notation.  When needed, this
dependence will be specified as follows: $I_n[f]$, $R_n[f]$, $I^\d_n[f]$,...

\subsubsection{Parameter partition}

Part of our work is to transfer information from the phase space of some
map $f \in \FF$ to a neighborhood of $f$ in the parameter space.
This is done in the following way.  We consider the
first landing map $L_i$: the complement of the domain of $L_i$ is a
hyperbolic Cantor set $K_i=I_i \setminus \cup C^\d_i$.
This Cantor set persists in a small parameter
neighborhood $J_i$ of $f$, changing in a continuous way.
Thus, loosely speaking, the domain of $L_i$ induces a persistent partition
of the interval $I_i$.

Along $J_i$, the first landing map is topologically the same (in a way
that will be clear soon).  However the critical value
$R_i[g](0)$ moves relative
to the partition (when $g$ moves in $J_i$).  This allows us to partition
the parameter piece $J_i$ in smaller pieces, each
corresponding to a region where $R_i(0)$ belongs to some fixed component of
the domain of the first landing map.

\begin{thm} [Topological Phase-Parameter relation]

Let $f \in \FF_\kappa$.
There is a sequence $\{J_i\}_{i \in \N}$ of nested
parameter intervals (the {\it principal parapuzzle nest} of $f$)
with the following properties.

\begin{enumerate}

\item $J_i$ is the maximal interval containing $f$ such that for all $g \in
J_i$ the interval
$I_{i+1}[g]=T^{(\kappa)}_{i+1}[g]$ is defined and changes in a
continuous way.  (Since the first return map to $R_i[g]$ has a central
domain, the landing map
$L_i[g]:\cup C^\d_i[g] \to I_{i+1}[g]$ is defined.)

\item $L_i[g]$ is topologically the same along $J_i$: there exists
homeomorphisms $H_i[g]:I_i \to I_i[g]$, such that
$H_i[g](C^\d_i)=C^\d_i[g]$.
The maps $H_i[g]$ may be chosen to change continuously.

\item There exists a homeomorphism $\Xi_i:I_i \to J_i$ such that
$\Xi_i(C^\d_i)$ is
the set of $g$ such that $R_i[g](0)$ belongs to $C^\d_i[g]$.

\end{enumerate}

\end{thm}

The homeomorphisms $H_i$ and $\Xi_i$ are not uniquely defined, it is easy to
see that we can modify them inside each $C^\d_i$ window keeping the above
properties.  However, $H_i$ and $\Xi_i$
are well defined maps if restricted to $K_i$.

This fairly standard phase-parameter result can be proved in many different
ways.  The most elementary proof is
probably to use the monotonicity of the
quadratic family to deduce the Topological Phase-Parameter relation from
Milnor-Thurston's kneading theory by purely combinatorial arguments.
Another approach is to use Douady-Hubbard's description of the combinatorics
of the Mandelbrot set (restricted to the real line) as does Lyubich
in \cite {parapuzzle} (see also \cite {AM4} for a more general case).

With this result we can define for any $f \in \FF_\kappa$ intervals
$J^j_i=\Xi_i(I^j_i)$ and $J^\d_i=\Xi_i(I^\d_i)$.  From the description we
gave it immediately follows that two intervals $J_{i_1}[f]$ and $J_{i_2}[g]$
associated to maps $f$ and $g$ are either disjoint or nested,
and the same happens for
intervals $J^j_i$ or $J^\d_i$.  Notice that if $g \in \Xi_i(C^\d_i) \cap
\FF_\kappa$ then $\Xi_i(C^\d_i)=J_{i+1}[g]$.

\subsubsection{Phase-Parameter relation}

In order to describe the metric properties of the phase-parameter map $\Xi$,
we will restrict ourselves to a smaller class of maps then $\FF_\kappa$, for
which we will be able to give a better description.  Those are maps for
which only finitely many returns $R_n$ are central, and are called {\it
simple maps} in \cite {AM1}.  We are able to restrict ourselves to this
class of maps due to the following result of Lyubich \cite {parapuzzle}:

\begin{thm}

Almost every map in $\FF$ has only finitely many
central returns in the principal nest.

\end{thm}

Even for simple maps, however, the regularity of $\Xi_i$ is not great:
there is too much dynamical information contained in it.  A solution to this
problem is to consider restrictions of $\Xi$ that ``forget''
some dynamical information.

\subsubsection{Geometric interpretation}

Before getting into those technical
details, it will be convenient to make an informal geometric
description of the topological statement we just made
and discuss in this context the difficulties
that will show up to obtain metric estimates.

The sequence of intervals $J_i$ is defined as the maximal parameter interval
containing $f$ satisfying two properties: the dynamical
interval $I_i$ has a continuation (recall that the boundary of
$I_i$ is preperiodic, so the meaning of continuation is quite clear),
and the first return map
to this continuation has always the same combinatorics.  Since it has the
same combinatorics, the partition $C^\d_i$ also has a continuation along
$J_i$.

Let us represent in two dimensions those continuations.  Let
$\II_i=\cup_{g \in J_i} \{g\} \times I_i[g]$ represent the
``moving phase space'' of
$R_i$.  It is a topological rectangle, its boundary consists of four
analytic curves, the top and bottom (continuations of the boundary points of
$I_i$) and the laterals (the limits of the continuations of $I_i$ as the
parameter converges to the boundary of $J_i$).
Similarly, the continuations of
each interval $C^\d_i$
form a strip $\CC^\d_i$ inside $\II_i$.  The resulting decomposition of
$\II_i$ looks like a flag with countable many strips.  The top and bottom
boundaries of those strips (and the strips themselves) are
horizontal in the sense that
they connect one lateral of $\II_i$ to the other.

\begin{rem}

The boundaries of the strips are more formally described as
forming a lamination in the topological rectangle $\II_i$, whose leaves are
codimension-one, and indeed
real analytic graphs over the first coordinate.  We remark that
in the complex setting,
the theory of codimension-one laminations is the same as the theory of
holomorphic motions, described in \cite {parapuzzle}, and which
is the basis of the actual phase-parameter analysis.

\end{rem}

Let us now look at the verticals $\{g\} \times I_i[g]$.  They are
all transversal to the strips of the flag, so we can consider the
``horizontal'' holonomy map between any two such verticals.  If we fix one
vertical as the phase space of $f$ while we vary the other,
the resulting family of
holonomy maps is exactly $H_i[g]$ as defined above.

Consider now the motion of $R_i[g](0)$ (the critical value
of the first return map to the continuation of $I_i$) inside $J_i$, which we
can represent by its graph $\DD=\cup_{g \in J_i} \{g\} \times
\{R_i[g](0)\}$. It is a diagonal to $\II_i$ in the sense
that it connects a corner
of the rectangle $\II_i$ to the
opposite corner.  In other words, if we vary continuously the
quadratic map $g$ (inside a slightly bigger parameter window then $J_i$),
we see the window $J_i$ appear when $g^{v_i}(0)$
enters $I_i[g]$ from one side and disappear when $g^{v_i}(0)$ escapes from
the other side (where $v_i$ is such that $R_i|_{I^0_i}=f^{v_i}$).

The main content of the Topological Phase-Parameter relation is
that the motion of the critical value is not
only a diagonal to $\II_i$ but to the
flag: it cuts each strip exactly once in a monotonic way with respect to the
partition.  Thus, the diagonal motion of the critical point is transverse in
a certain sense to the horizontal motion of the partition of the phase
space (strips).  The phase-parameter map is just the
composition of two maps:
the holonomy map between two transversals to the flag
(from the ``vertical'' phase space of $f$ to the diagonal $\DD$)
followed by
projection on the first coordinate (from $\DD$ to $J_i$).

\begin{rem}

One big advantage of complex analysis is that ``transversality can be
detected for topological reasons''.  So, while the statement that the
critical point goes from the bottom to the top of $\II_i$ does not imply
that it is transverse to all horizontal strips, the corresponding
implication holds for the complex analogous of those statements.  This is a
consequence of the Argument Principle.

\end{rem}

Let us now pay attention to the geometric format of those strips.
The set $\cup_{g \in J_i} \{g\} \times
R_i[g](I^0_i[g])$ is a topological triangle formed by the diagonal $\DD$,
one of the laterals of $\II_i$ (which we will call
the right lateral\footnote {It
is possible to prove that it is indeed located at the right side (with the
usual ordering of the real line).})
and either the top or bottom of $\II_i$.  In particular, the strip $\II^0_i$
is not a rectangle, but a triangle: the left side of $\II^0_i$ degenerates
into a point.  By their dynamical definition, all strips $C^\d_i$ also
share the same property: the $C^\d_i[g]$ are collapsing as $g$ converges to
the left boundary of $J_i$.  In particular the partition of the phase space
of $f$ must be metrically very different from the partition of the phase
space of some $g$ close to the left boundary of $J_i$.

This shows that it is not reasonable to expect the phase-parameter map
$\Xi_i|_{K_i}$ to be very regular (uniformly H\"older for
instance\footnote{The Hausdorff dimension of $K_i[g]$ is not constant for $g
\in J_i$, so Lipschitz estimates are certainly out of reach.}): if it
was true that the phase-parameter relation is always regular,
then the phase partitions of $f$ and $g$ would have to be metrically
similar (since a correspondence between both partitions can be obtained
as composition the phase-parameter
relation for $f$ and the inverse of the phase-parameter relation for $g$).

Let us now consider the decomposition of $\II_i$ in strips $\II^j_i$ (the
continuations of $I^j_i$).  This new flag is rougher than the previous one:
each of its strips $\II^j_i$, $j \neq 0$ can be obtained as
the (closure of the) union of $\CC^\d_i$ where $\d$ starts with $j$. 
However, the strips are nicer: they are indeed rectangles if $j \neq 0$,
though the ``niceness'' gets weaker and weaker as we get closer to the
central strip.  This suggests one way to obtain a regular map from $\Xi_i$:
work with the rougher partition $I^j_i$ outside of a certain small
neighborhood of the critical strip (this neighborhood will be
introduced in the next section, it will be called the gape interval).
This procedure will indeed have the desired effect in the sense that we will
be able to prove that for simple maps $f$
the restriction of $\Xi_i$ to $I_i
\setminus \cup I^j_i$ has good regularity outside of the gape interval
(this is PhPa2 in the Phase-Parameter relation below).

The resulting estimate does not say anything about what happens inside the
rough partition by $J^j_i$.  To do so, we consider the finer flag (whose
strips are the $\CC^\d_i$) intersected with
the rectangles $\QQ^j_i=\cup_{g \in J^j_i} \{g\} \times I^j_i[g]$ (those
rectangles cover the diagonal $\DD$
formed by the motion of the critical value).
While the strips degenerate near the left boundary point of $J_i$, they
intersect each $\QQ^j_i$ in a nice rectangle (or the empty set).
It will be indeed possible to prove that
the phase-parameter map restricted to those rectangles
is quite regular in the sense that if $f$ is a simple map such that
$f \in J^j_i$ (that is, $(f,R_n(0)) \in \QQ^j_i$),
the restriction of $\Xi_i$ to $I^j_i$ has good
regularity (this is PhPa1 in the Phase-Parameter relation below).

\subsection{Quasisymmetric maps}

As we just described,
phase-parameter maps can be viewed as holonomy maps of
``flags'', which are codimension-one laminations with real analytic leaves. 
It turns out that such objects inherit some ``automatic'' regularity
from their complexifications: they are quasisymmetric, at least away from
the boundary (where we have the bad effects we just described).  The theory
of quasisymmetric maps is a well developed subject, but we will need just
the definition and a couple of elementary properties.

Let $k \geq 1$ be given.
We say that a homeomorphism $f:\R \to \R$ is {\it quasisymmetric}
with constant $k$ if for all $h > 0$
$$
\frac {1} {k} \leq \frac {f(x+h)-f(x)} {f(x)-f(x-h)} \leq k.
$$

The space of quasisymmetric maps is a group under composition, and the set
of quasisymmetric maps with constant $k$ preserving a given interval is
compact in the uniform topology of compact subsets of $\R$.  It also follows
that quasisymmetric maps are H\"older.  Quasisymmetric maps are much better
than H\"older though: the key additional property, used to no end in the
statistical analysis is that the definition of quasisymmetric maps (and
associated constants) is {\it scaling invariant} (invariant under affine
changes of coordinates.

To describe further the properties of quasisymmetric maps, we need the
concept of quasiconformal maps and dilatation
so we just mention a result of Ahlfors-Beurling
which connects both concepts: any quasisymmetric map extends
to a quasiconformal real-symmetric map of
$\C$ and, conversely, the
restriction of a quasiconformal real-symmetric map of $\C$ to $\R$ is
quasisymmetric.  Furthermore, it is possible to work out upper bounds
on the dilatation $\g$ (of an optimal extension) depending only on $k$ and
conversely: it turns out that
$\g$ is close to $1$ if and only if $k$ is close to $1$.

The constant $k$ is awkward to work with: the inverse of a quasisymmetric
map with constant $k$ may have a larger constant.  We will therefore work
with a less standard constant: we will say that $h$ is
$\g$-quasisymmetric ($\g$-qs) if $h$ admits a quasiconformal symmetric
extension to $\C$ with dilatation bounded by $\g$.  This definition
behaves much better: if $h_1$ is $\g_1$-qs and $h_2$ is
$\g_2$-qs then $h_2 \circ h_1$ is $\g_2 \g_1$-qs.

If $X \subset \R$ and $h:X \to \R$ has a $\g$-quasisymmetric extension to
$\R$ we will also say that $h$ is $\g$-qs.

\subsubsection{The Phase-Parameter relation}

As we discussed before, the dynamical information contained in $\Xi_i$
is entirely given by $\Xi_i|_{K_i}$: a map obtained by $\Xi_i$ by
modification inside a $C^\d_i$ window has still the same properties. 
Therefore it makes sense to ask about the regularity of $\Xi_i|_{K_i}$.  As
we anticipated before we must erase some information to obtain good results. 

If $i>1$, we define the {\it gape interval} $\tilde I_{i+1}$ as follows.
Let $\d$ be such that
$R_i|_{I_{i+1}}=L_{i-1} \circ R_{i-1}=R^\d_{i-1} \circ R_{i-1}$, so that
$I_{i+1}=(R_{i-1}|_{I_i})^{-1}(C^{\d}_{i-1})$.  The gape interval is defined
as $\tilde I_{i+1}=(R_{i-1}|_{I_i})^{-1}(I^\d_{i-1})$.
Notice that $I_{i+1} \subset \tilde I_{i+1} \subset I_i$.  Furthermore,
for each $I^j_i$, the gape interval
$\tilde I_{i+1}$ either contains or is disjoint from $I^j_i$.

Let $f \in \FF_\kappa$ and let $\tau_i$ be such that
$R_i(0) \in I^{\tau_i}_i$.  We define two Cantor sets,
$K^\tau_i=K_i \cap I^{\tau_i}_i$ which
contains refined information restricted to the $I^{\tau_i}_i$ window and
$\tilde K_i=I_i \setminus (\cup I^j_i \cup \tilde I_{i+1})$, which
contains global information, at the cost of erasing information inside each
$I^j_i$ window and in $\tilde I_{i+1}$.

\begin{thm} [Phase-Parameter relation]

Let $f$ be a simple map.  For all $\g>1$ there exists $i_0$ such that
for all $i>i_0$ we have

\begin{description}

\item[PhPa1] $\Xi_i|_{K^\tau_i}$ is $\g$-qs,

\item[PhPa2] $\Xi_i|_{\tilde K_i}$ is $\g$-qs,

\item[PhPh1] $H_i[g]|_{K_i}$ is $\g$-qs if $g \in J^{\tau_i}_i$,

\item[PhPh2] the map $H_i[g]|_{\tilde K_i}$ is
$\g$-qs if $g \in J_i$.

\end{description}

\end{thm}

The proof of the Phase-Parameter relation is based on complex methods, and
the ideas involved go beyond the scope of this note.
The ideas which are necessary in the analysis come from the work of
Lyubich in \cite {parapuzzle}, where a general method based on the
theory of holomorphic motions was introduced to deal with this kind of
problem.  A sketch of the derivation
of the specific statement of the Phase-Parameter relation from the general
method of Lyubich was given in the
Appendix A of \cite {AM1}.  The reader can find full details
(in a more general context
than quadratic maps) in \cite {AM4}.


\subsection{Basic ideas of the statistical analysis}

We will now describe the statistical analysis done in \cite {AM1}.  One of
the key difficults to overcome is the fact that the Phase-Parameter relation
is not Lipschitz (quasisymmetric maps are not even absolutely continuous in
general).  Thus, instead of working with Lebesgue measure in phase space, we
are lead to work with ``quasisymmetric capacities'' defined as follows: if
$\g>1$, the $\g$-qs capacity of a set $X$ in an interval $T$ is
\be
p_\g(X|T)=\sup \frac {|h(X)|} {|h(T)|}
\ee
where the suppremum is taken over all $\g$-qs maps $h:\R \to \R$.

By design, sets of small capacity in phase space must be taken by the
phase-parameter map to sets of small Lebesgue measure in parameter space.
There is a price to be paid: capacities are not probabilities (one may have
two disjoint sets with capacities close to $1$), so we must do some work in
order to be able to apply statistical laws as the Law of Large Numbers
and the Law of Large Deviations.

The fact that $\g$-quasisymmetric maps are H\"older (with good constants if
$\g$ is close to $1$) is not quite enough to do any analysis: scaling
invariance is an important part of the renormalization game.  We will
actually exploit scaling invariance through
the following property of capacities: if $T^j \subset T$ is a
disjoint family of intervals covering $T \cap X$, then
\be
p_\g(X|T) \leq p_\g(\cup T^j|T) \sup_j p_\g(X|T^j),
\ee
which fits particularly well with the tree structure of the family
$\{I^\d_n\}_{\d \in \Omega}$ (organized by inclusion).

\subsubsection{Borel-Cantelli and the parameter exclusion process}

The parameter exclusion process consists in obtaining successively
smaller (but still full-measure) classes of maps for which we can give
a progressively refined statistical description of the dynamics.
This is done inductively as follows: we pick a class $X$ of maps
(which we have previously shown to
have full measure among non-regular maps) and for each map in $X$
we proceed to describe the dynamics (focusing on the statistical behavior
of return and landing maps for deep levels of the principal nest),
then we use this information to show that a subset $Y$ of $X$
(corresponding to parameters for which the statistical behavior of the
{\it critical orbit} is not anomalous) still has full measure.
An example of this parameter exclusion process is done by Lyubich in \cite
{parapuzzle} where he shows using a probabilistic argument that
the class of simple maps has full measure in $\FF$.

Let us now describe our usual argument (based on the argument of Lyubich
which in turn is a variation of the Borel-Cantelli Lemma).
Assume that at some point we know how to prove that almost every
simple map belongs to a certain set $X$.
Let $Q_n$ be a (bad) property that a map
may have (usually some anomalous statistical parameter
related to the $n$-th stage of the principle nest).
Suppose we prove that if $f \in X$ then
the probability that a map in $J_n(f)$ has the property $Q_n$ is bounded by
$q_n(f)$ which is shown to be summable
for all $f \in X$.  We then conclude that almost every map does not have
property $Q_n$ for $n$ big enough.

Sometimes we also apply the same argument, proving instead that $q_n(f)$
is summable where $q_n(f)$ is the
probability that a map in $J^{\tau_n}_n(f)$ has property $Q_n$,
(recall that $\tau_n$ is such that $f \in J^{\tau_n}_n(f)$).

In other words, we apply the following simple general result.

\begin{lemma} \label {measureth}

Let $X \subset \R$ be a measurable set such that for each
$x \in X$ is defined a sequence
$D_n(x)$ of nested intervals converging to $x$
such that for all $x_1,x_2 \in X$
and any $n$, $D_n(x_1)$ is either equal or disjoint to $D_n(x_2)$.  Let
$Q_n$ be measurable subsets of
$\R$ and $q_n(x)=|Q_n \cap D_n(x)|/|D_n(x)|$.  Let $Y$ be
the set of all $x \in X$ which belong to at most
finitely many $Q_n$.
If $\sum q_n(x)$ is finite for almost any $x \in X$ then $|Y|=|X|$.

\end{lemma}

In practice, we will estimate the capacity of sets in the phase space:
that is, given a map $f$ we will obtain subsets $\tilde Q_n[f]$ in the phase
space, corresponding to bad branches of return or landing maps.  We will
then show that for some $\g>1$ we have
$\sum p_\g(\tilde Q_n[f]|I_n[f])<\infty$ or
$\sum p_\g(\tilde Q_n[f]|I^{\tau_n}_n[f])<\infty$.  We will then use PhPa2
or PhPa1, and the measure-theoretical lemma above to conclude that with
total probability among non-regular maps, for all $n$ sufficiently big,
$R_n(0)$ does not belong to a bad set.

\subsubsection{A case study}

We will now describe in detail how to apply the measure-theoretical
argument and the Phase-Parameter relation.  In order to illustrate our
ideas, we will discuss informally the first statistical result of \cite
{AM1}, which is quite simple yet particularly important for our strategy.

For a map $f \in \Delta_\kappa$ (recall that, as always,
we work in a fixed level $\kappa$ of renormalization),
let us associate a sequence of ``statistical parameters'' in some way.
A good example of statistical parameter is $s_n$, which denotes the
number of times the critical point $0$ returns to $I_n$ before the first
return to $I_{n+1}$.  Each of the points of the
sequence $R_n(0)$,...,$R_n^{s_n}(0)$ can be located anywhere
inside $I_n$.  Pretending that the distribution of those points is indeed
independent and uniform with respect to Lebesgue measure, we may expect
that typical values of $s_n$ concentrate near (in an appropriate sense)
$c_n^{-1}$, where $c_n=|I_{n+1}|/|I_n|$.  For the ``random model'',
{\it near} may be interpreted in logharithmic scale in terms of the
difference
\be
\frac {\ln s_n} {\ln c_n^{-1}}-1,
\ee
and one sees indeed that the concentration is more marked the smaller
$c_n$ is (of course for the random model one can obtain much better
than loharithmic estimates).

Let us try to make such an estimate rigorous.  Consider the set of
points $A_k \subset I_n$ which iterate exactly $k$ times in $I_n$ before
entering $I_{n+1}$.  Then most points $x \in I_n$ belong to some
$A_k$ with $k$ in a (logharithmic) neighborhood of $c_n^{-1}$ (if we forget
about distortion, the probability of $A_k$ is $c_n (1-c_n)^k$).
By most, we mean that the complementary event has small
probability, say $q_n$, for some summable sequence $q_n$.
This neighborhood has to be computed precisely using a statistical argument.
In this case, if we choose the neighborhood
$c_n^{-1+2\epsilon}<k<c_n^{-1-\epsilon}$, we obtain the sequence
$q_n<c_n^\epsilon$ which is indeed summable for all simple maps
$f$ by \cite {attractors}.

If the phase-parameter relation was Lipschitz, we would now argue as
follows: the probability of a parameter be such that $R_n(0) \in A_k$ with
$k$ out of the ``good neighborhood'' of values of $k$ is also summable
(since we only multiply those probabilities by the Lipschitz constant) and
so, by the measure-theoretical argument of Lemma \ref {measureth},
for almost every parameter this only happens a finite number of times.

Unfortunately, the Phase-Parameter relation is not Lipschitz.  To make the
above argument work, we must have better control of the size of the ``bad
set'' of points which we want the critical value $R_n(0)$
to not fall into.  In order to do so, in the statistical analysis of the
sets $A_k$, we control instead the quasisymmetric capacity of the
complement of points falling in the good neighborhood.  This makes
the analysis sometimes much more difficult since capacities are not
probabilities. 
This will usually introduce some error that was not present in the naive
analysis, leading to the $\epsilon$ in the range of
exponents present above.
This is why we do not try to do better than
estimates in logharithmic scale: if we were
not forced to deal with capacities, we could get much finer estimates.

Incidentally, to keep the error low, making $\epsilon$ close to $0$,
we need to use capacities with constant $\g$ close to $1$.
Fortunately, our
Phase-Parameter relation has a constant converging to $1$, which will allow
us to partially get rid of this error.
(Indeed, with $\g$
close to $1$, we can get $p_\g(A_k|I_n)<c_n^{1-\delta}(1-c_n^{1+\delta})^k$
with $\delta$ close to $0$, which will be enough for our purposes.)

Coming back to our problem, we see that we should concentrate in proving
that for almost every parameter, certain $\epsilon$-bad sets have summable
$\g$-qs capacities for some constant $\g$ independent of $n$ (but which can
depend on $f$ and $\epsilon$).

There is one final detail we should pay attention to: there are two
phase-parameter statements, and we should use the right one.  More
precisely, there will be situations where we are analyzing some sets which
are union of $I^j_n$ (return sets), and sometimes the relevant sets are
union of $C^\d_n$ (landing sets).  In the first case, we should use the
PhPa2 and in the second the PhPa1.  Notice that
our phase-parameter estimates only allow us to ``move the
critical point'' inside $I_n$ with respect to the partition by $I^j_n$: to
do the same with respect to the partition by $C^\d_n$, we must restrict
ourselves to $I^{\tau_n}_n$.
In all cases, however, the bad sets considered
should be either union of $I^j_n$ or $C^\d_n$.

For our specific example, since the $A_k$ are union of $C^\d_n$, we must use
PhPa1.  In particular we have to study the capacity of a bad set inside
$I^{\tau_n}_n$.  Here is the estimate that we should go after (see Lemma 4.2
of \cite {AM1} for a more precise statement):

\begin{lemma}

For almost every parameter,
for every $\epsilon>0$, there exists $\g$ such that
$p_\g(X_n|I^{\tau_n}_n)$
is summable, where $X_n$ is the set of points $x \in I_n$ which enter
$I_{n+1}$ either before $c_n^{-1+\epsilon}$ or after $c_n^{-1-\epsilon}$
returns to $I_n$.

\end{lemma}

We are now in position to use PhPa1 to make the corresponding
parameter estimate:
using the measure-theoretic argument,
we get (in Lemma 4.3 of \cite {AM1}) that with total probability
\be \label {s_n}
\lim_{n \to \infty} \frac {\ln s_n} {\ln c_n^{-1}}=1.
\ee

This particular estimate we chose to describe in this section is extremely
important for the analysis to follow:
we can use $s_n$ to estimate $c_{n+1}$ directly from below:
\be
\frac {\ln c_{n+1}^{-1}} {s_n} \to \infty
\ee
so this last lemma implies (Corollary 4.4 of \cite {AM1}) that
$c_n^{-1}$ grows at least as fast as a tower of $2$'s of height $n-C$ for
somes $C>0$ independent of $n$
(this kind of decay/growth will be called torrential).

For general simple maps,
the best information is given by \cite {attractors}:
$c_n$ decays exponentially
(this was actually used to obtain summability of $q_n$
in the above argument).  This improvement from
exponential to torrential should give the reader an idea of
the power of this kind of statistical analysis.

\subsection{Collet-Eckmann and polynomial recurrence: strategy of the proof}

We now describe the key ideas involved in the proof of the main results of
\cite {AM1}, namely: Almost every non-regular quadratic map $f$ satisfies
the Collet-Eckmann condition
\be
\liminf_{n \to \infty} \frac {\ln|Df^n(f(0))|} {n}>0
\ee
and the orbit of the critical point has
polynomial recurrence with exponent 1, that is
\be
\limsup_{n \to \infty} \frac {-\ln |f^n(0)|} {\ln n}=1.
\ee

\subsubsection{Distribution of the hyperbolicity random variable}

Let us first explain how the information on statistical
parameters of a typical non-regular map $f$ can be used to obtain estimates
of hyperbolicity along the critical orbit that imply the Collet-Eckmann
condition.  We make several simplifications, in particular we
don't discuss here
the difficulty involved in working with capacities instead of probabilities.

Let us start by thinking of hyperbolicity at a given level $n$ as a
random variable $\lambda_n(j)$ (introduced in \S 7.2 of \cite {AM1})
which associates to each non-central branch of
$R_n$ its average expansion, that is, if $R_n|_{I^j_n}=f^{r_n(j)}$ ($r_n(j)$
is the return time of $R_n|_{I^j_n}$), we let
$\lambda_n(j)=\ln |Df^{r_n(j)}|/r_n(j)$ evaluated at some point
$x \in I^j_n$, say, the point where $|Df^{r_n(j)}|$ is minimal\footnote{
The choice of the point in $I^j_n$ turns out to be not very relevant
because it is possible to obtain reasonable (polynomial)
``almost sure'' bounds on distortion, see Lemma 4.10 of \cite {AM1}.}.

Our tactic is to evaluate the evolution of the distribution of
$\lambda_n(j)$ as $n$ grows.  The basic information we will use to start our
analysis is the hyperbolicity estimate of Lemma \ref {hyperbol},
which, together with our
distortion estimates, shows that $\lambda_n=\inf_j \lambda_n(j)>0$
for $n$ big enough.  We then fix such a big level $n_0$ and
the remaining of the analysis will be based on inductive
statistical estimates for levels $n>n_0$.

Of course, nothing guarantees a priori that $\lambda_n$ does not
decay to $0$.  Indeed, it turns out that $\liminf_{n \to \infty}
\lambda_n>0$, but as a consequence of
the Collet-Eckmann condition and our distortion estimates.
But this is not what we will
analyze: we will concentrate on showing that $\lambda_n(j)>\frac {n+1} {2n}
\lambda_{n_0}>\frac {1} {2} \lambda_{n_0}$
outside of a ``bad set'' of torrentially small $\g$-qs capacity.  The
complementary set of hyperbolic branches will be called good.

To do so, we inductively describe branches of level $n+1$ as compositions of
branches of level $n$.  Assuming that most branches of level $n$ are good,
we consider branches of level $n+1$ which spend most of their time in good
branches of level $n$.  They inherit hyperbolicity from good branches of
level $n$, so they are themselves good of level $n+1$.
To make this idea work we should also have additionally a condition of
``not too close returns'' to avoid drastic reduction of
derivative due to the critical point.

The fact that most branches of level $n$ were good (quantitatively: branches
which are not good have capacity bounded by some small $q_n$) should
reflect on the fact that most branches of level $n+1$ spend a small
proportion of their time (less than $6 q_n$) on branches which
are not good, and so most branches of level $n+1$ are also good
(capacity of the complement is a small $q_{n+1}$): indeed
the notion of most should improve from level to level, so that
$q_{n+1} \ll q_n$ (in order for this argument to work, the hyperbolicity
requirements in the notion of good must become slightly more flexible when
we go from level to level).
This reflects the tendency of averages of random variables to
concentrate around the expected
value with exponentially small errors (Law of Large numbers and Law of Large
deviations).  Those laws give better results if we average over a larger
number of random variables.  In particular,
those statistical laws are very effective in our case,
since the number of random variables that we average
will be torrential in $n$: our arguments will typically lead to estimates as
$\ln q_{n+1}^{-1}>q_n^{-1+\epsilon}$ (torrential decay of $q_n$).

In practice, we will obtain good branches in a more systematic way.  We
extract from the above crude arguments a couple of features that should
allow us to show that some branch is good.  Those features define what we
call a very good branch:
\begin{enumerate}
\item for very good branches we can control
the distance of the branch to $0$ (to avoid drastic loss of derivative);
\item the definition of very good branches has
an inductive component: it must be a composition of many branches,
most of which are themselves very good of the previous level
(with the hope of propagating hyperbolicity inductively);
\item the distribution of return times of branches of the previous level
taking part in a very good branch has a controlled ``concentration around
the average''.
\end{enumerate}

Let us explain the third item above: to compute the hyperbolicity of a
branch $j$ of level $n+1$, which is a composition of several branches
$j_1,...,j_m$ of level $n$ we are essentially estimating
$$
\frac {\sum r_n(j_i) \lambda_n(j_i)} {\sum r_n(j_i)},
$$
the ratio between the total expansion and the total time of the branch.
The second item assures us that many branches $j_i$ are very good, but this
does not mean that their total time is a reasonable part of
the total time $r_{n+1}(j)$ of the branch $j$.  This only holds if
we can guarantee some
concentration (in distribution) of the values of $r_n(j)$.

With those definitions we can prove that very good branches are good, but to
show that very good branches are ``most branches'', we need to understand
the distribution of the return time random variable
$r_n(j)$ that we discuss later.

Let us remark that our statistical work so far (showing that
good (hyperbolic) branches are most branches) is not yet enough to conclude
Collet-Eckmann: indeed we have controlled
hyperbolicity only at full returns.  To estimate hyperbolicity at any moment
of some orbit, we must use good branches as building blocks of hyperbolicity
of some special branches of landing maps (cool landings).  Branches which
are not very good are sparse inside truncated cool landings, so that if we
follow a piece of orbit of a point inside a cool landing (not necessarily up
to the end), we still have enough hyperbolic blocks to estimate the growth
of derivative.

After all those estimates, we use the Phase-Parameter relation to move the
critical value into cool landings, and obtain exponential growth of
derivative of the critical value (with rate bounded from below by
$\lambda_{n_0}/2$).

\subsubsection{The return time random variable}

As remarked above, to study the hyperbolicity random variable
$\lambda_n(j)$, we must first estimate
the distribution of the return time
random variable $r_n(j)$.  It is worth to discuss some key ideas of this
analysis (\S 6 of \cite {AM1}).

Intuitively, the ``expectation'' of $r_n(j)$ should be concentrated in a
neighborhood of $c_{n-1}^{-1}$:
pretending that iterates $f^k(x)$ are
random points in $I$, we expect to wait about $|I|/|I_n|$ to get back to
$I_n$.  But $|I_n|/|I|=c_{n-1}c_{n-2}...c_1 |I_1|/|I|$.  Since $c_n$ decays
torrentially, we can estimate $|I_n|/|I|$ as $c_{n-1}^{1-\epsilon}$ (it is
not worth to be more precise, since $\epsilon$ errors in the exponent will
appear necessarily when considering capacities).  Although this naive
estimate turns out to be true, we of course don't try to follow this
argument: we never try to iterate $f$ itself, only return branches.

The basic information we use to start is again the hyperbolicity estimate of
Lemma \ref {hyperbol}.  This information gives us exponential tails
for the distribution of $r_n(j)$ (the $\g$-qs capacity of
$\{r_n(j)>k\}$ decays exponentially in
$k$).  Of course we have no information on the exponential rate: to control
it we must again use an inductive argument which studies the propagation of
the distribution of $r_n$ from level to level.  The idea again
is that random variables add well and the relation between $r_n$ and
$r_{n+1}$ is additive: if the branch $j$ of level $n+1$ is the composition
of branches $j_i$ of level $n$, then
$r_{n+1}(j)$ is the sum of the $r_n(j_i)$.

Using that the transition from level to level involves
adding a large number of random variables (torrential), we are able to
give reasonable bounds for the decay of the tail of $r_n(j)$ for $n$ big
(this step is what we call a Large Deviation estimate).  Once we control
this tail, an estimate of the concentration of the
distribution of return times becomes natural
from the point of view of the Law of Large Numbers.

\subsubsection{Recurrence of the critical orbit}

After doing the preliminary work on the distribution of return times,
the idea of the estimate on recurrence (which is done
in \S 8.2 of \cite {AM1}) is quite transparent.

We first estimate the rate that a typical sequence $R_n^k(x)$ in $I_n$
approaches $0$, before falling into $I_{n+1}$.  If the sequence $R_n^k(x)$
was random, then this recurrence would clearly be polynomial with exponent
$1$.  The system of non-central
branches is Markov with good estimates of distortion, so it is no surprise
that $R_n^k(x)$ has the same recurrence properties, even if the system is
not really random.  We can then conclude that some inequality as
\begin{equation} \label {bound rec}
|R_n^k(x)|>|I_n|2^{-n}k^{-1-\epsilon}
\end{equation}
holds for most orbits (summable complement).

We must then relate the recurrence in terms of iterates of $R_n$ to the
recurrence in terms of iterates of $f$.  Since in the Collet-Eckmann
analysis we proved that (almost surely) the
critical value belongs to a cool landing, it is enough to do the estimates
inside a cool landing.  But cool landings are formed by well distributed
building blocks with good distribution of return times, so we can relate
easily those two recurrence estimates.

To see that when we pass from the estimates in terms of iterations by
$R_n$ to iterations in term of $f$ we still get polynomial recurrence,
let us make a rough estimate which indicates that
$R_n(0)=f^{v_n}(0)$ is at distance approximately $v_n^{-1}$ of $0$. 
Indeed $R_n(0)$ is inside
$I_n$ by definition, so we have a trivial upper bound
$|R_n(0)|<c_{n-1}$.  Using the
phase-parameter relation, the critical orbit has controlled recurrence (in
terms of (\ref {bound rec})), thus we get
$|R_n(0)|>2^{-n}|I_n|>c_{n-1}^{1+\epsilon}$.  Together with the upper bound,
this implies that $|R_n(0)|$ is of order
$c_{n-1}^{-1}$.  On the
other hand, $v_n$ (number of iterates of $f$ before getting to $I_n$)
is at least $s_{n-1}$ (number of iterates of $R_{n-1}$ before getting to
$I_n$).  According to (\ref {s_n}), $s_{n-1}$ is of order $c_{n-1}^{-1}$,
so this
argument gives the lower bound $v_n>c_{n-1}^{-1+\epsilon}$.  On the other
hand, $v_n$ is $s_{n-1}$ times the average time of branches $R_{n-1}$: due
to our estimates on the distribution of return times,
$$
v_n<s_{n-1} c_{n-2}^{1-\epsilon}<c_{n-1}^{-1-\epsilon}c_{n-2}^{-1-\epsilon}<
c_{n-1}^{-1-2\epsilon}
$$
(here we use that $0$ is a ``typical'' point for the distribution of return
times, since it falls in cool landings).  Together with the lower bound,
this implies that $v_n$ is of order $c_{n-1}^{-1}$, and we get
$$
1-4\epsilon<\frac {\ln |R_n(0)|} {\ln v_n}<1+4\epsilon.
$$

\subsubsection{Some technical details}

The statistical analysis described above is considerably
complicated by the use of capacities: while traditional
results of probability can be used as an
inspiration for the proof (as outlined here),
we can not actually use them.
We also have to use statistical
arguments which are adapted to tree decomposition of landings into returns:
in particular, more sophisticated analytic estimates
are substituted by more ``bare-hands'' techniques.

Following the details of the actual proof in \cite {AM1},
the reader will notice that we work very often with a sequence of
quasisymmetric constants which decrease
from level to level but stays bounded away from $1$.
We don't work with a fixed capacity because, when adding
random variables as above, some distortion is introduced.  We can make the
distortion small but not vanishing,
and the distortion affects the constant of the
next level: if we could make estimates of distribution using some constant
$\g_n$, in the next level the estimates are in terms
of a smaller constant $\g_{n+1}$.  These ideas are introduced
in \S 5 of \cite {AM1}.

Since the phase-parameter relation has two parts, our statistical analysis
of the transition between two levels will very often involve two steps: one
in order to move the critical value out of bad branches of the return map
$R_n$, and another to move it inside a given branch of
$R_n$ outside of bad branches of the landing map $L_n$.

Fighting against the technical difficulties is the torrential decay of
$c_n$.  The typical values of statistical parameters appearing in the
analysis of level $n$
are usually related to $c_n$ or $c_{n-1}$, up to a small error in the
exponent.  When statistical parameters of different levels interact
usually only one
of them will determine the order of magnitude of the
result.  This is specially true since all our
estimates include an $\epsilon$ error in the exponent.  The reader should
get used to estimates as ``$c_n c_{n-1}$ is approximately $c_n$'', in the
sense that the ratio of the
logarithms of both quantities is actually close to $1$ (compare the
estimates in the end of the last section, specially relating $s_n$ and
$v_n$).  Even if many proofs are quite technical,
they are also quite robust due to this.

\comm{
\subsection{Propagation of small errors}

\subsection{Renormalization versus inducing}

Let us finish this outline section with a comment on the philosophy of the
our method (Lyubich's generalized renormalization) compared to the method of
Jakobson (inducing).

While both methods attack, as expected, the critical region, the inducing
method can be roughly described as ``attempting to go back to large
scales''.  In other words, one looks.
}


\begin{thebibliography}{*****}


\bibitem[A]{A} A. Avila.  Bifurcations of unimodal maps: the topological and
metric picture.  Thesis IMPA (2001)
(www.math.sunysb.edu/$\sim$artur).

\bibitem[ALM]{ALM} A. Avila, M. Lyubich, W. de Melo.  Regular or
stochastic dynamics in real analytic families of unimodal maps.
Preprint (www.math.sunysb.edu/$\sim$artur).  To appear in Inventiones Math.

\bibitem[AM1]{AM1} A. Avila, C. G. Moreira.  Statistical properties of
unimodal maps: the quadratic family.  To appear in Annals of Math.

\bibitem[AM2]{AM2} A. Avila, C. G. Moreira.  Statistical properties of
unimodal maps: smooth families with negative Schwarzian derivative.
Preprint (www.arXiv.org).  To appear in Ast\'erisque.

\bibitem[AM3]{AM4} A. Avila, C. G. Moreira.  Phase-Parameter relation and
sharp statistical properties in general families of unimodal maps.
In preparation.

\bibitem[BBM]{BaBeM} V. Baladi, M. Benedicks, V. Maume.  Almost sure rates
of mixing for i.i.d. unimodal maps.  Preprint (1999), to appear Ann. E.N.S.

\bibitem[BV]{BV} V. Baladi, M. Viana.  Strong stochastic stability and
rate of mixing for unimodal maps.  Ann. scient. \'Ec. Norm. Sup., v. 29
(1996), 483-517.

\bibitem[BC1]{BC1} M. Benedicks, L. Carleson. On iterations of $1-ax^2$
 on (-1,1). Ann. Math., v. 122 (1985), 1-25. 

\bibitem[BC2]{BC2}  M. Benedicks, L. Carleson.
 On dynamics of the H\'enon map. Ann. Math., v. 133 (1991), 73-169. 


\bibitem[BLS]{BLS} H. Bruin, S. Luzzatto, S. van Strien.  Decay of
correlations in one-dimensional dynamics.  Preprint (www.arXiv.org).
To appear in Ann. Sci. ENS.


\bibitem[GS1]{GS} J. Graczyk, G. Swiatek.  Generic hyperbolicity in the
logistic family.  Ann. of Math., v. 146 (1997), 1-52.

\bibitem[GS2]{GS2} J. Graczyk, G. Swiatek.  Induced expansion for quadratic
polynomials.  Ann. Sci. Éc. Norm. Supér., IV. Sér. 29, No.4 (1996),
399-482. 

\bibitem[HK]{HK} F. Hofbauer, G. Keller. Quadratic maps without
asymptotic measure.  Comm. Math. Physics, v. 127 (1990), 319-337.

\bibitem[J]{J} M. Jacobson. Absolutely continuous invariant measures
for one-parameter families of one-dimensional maps. Comm. Math. Phys.,
 v. 81 (1981), 39-88. 

\bibitem[Jo]{Jo} S. D. Johnson.  Singular measures without restrictive
intervals.  Comm. Math. Phys., 110 (1987), 185-190.

\bibitem[KN]{KN} G. Keller, T. Nowicki.  Spectral theory, zeta functions and
the distribution of periodic points for Collet-Eckmann maps.  Comm. Math.
Phys., 149 (1992), 31-69.


\bibitem[L1]{attractors} M. Lyubich. Combinatorics, geometry and attractors of
     quasi-quadratic maps. Ann. Math, {\bf 140} (1994), 347-404.

\bibitem[L2]{puzzle} M. Lyubich. Dynamics of quadratic polynomials, I-II. 
   Acta Math., {\bf 178} (1997), 185-297. 

\bibitem[L3]{parapuzzle} M. Lyubich. Dynamics of quadratic polynomials, III.
Parapuzzle and SBR measure.
Asterisque,
v. 261 (2000),  173 - 200.

\bibitem[L4]{universe} M. Lyubich.  Feigenbaum-Coullet-Tresser
   universality and Milnor's hairiness conjecture. Ann. of Math. (2) 149
   (1999), no. 2, 319--420.

\bibitem[L5]{regular} M. Lyubich. Almost every real quadratic map
is either regular or stochastic. Ann. of Math. (2) 156 (2002), no. 1, 1-78.

\bibitem[L6]{cambridge} M. Lyubich. Renormalization ideas in conformal
   dynamics. Current developments in mathematics, 1995 (Cambridge, MA),
   155--190, Internat. Press, Cambridge, MA, 1994.

\bibitem[L7]{notices} M. Lyubich. The quadratic family as a
   qualitatively solvable model of chaos. Notices Amer. Math. Soc. 47
   (2000), no. 9, 1042--1052.

\bibitem[MN]{MN} M. Martens, T. Nowicki. Invariant measures for
Lebesgue typical quadratic maps. 
Asterisque,
v. 261 (2000), 239 - 252. 

\bibitem[MSS]{MSS} R. Ma\~n\'e, P. Sad \& D. Sullivan. 
 On the dynamics of rational maps, Ann. scient. Ec. Norm. Sup., 
  {\bf 16}  (1983), 193-217.

\bibitem[MvS]{MvS} W. de Melo, S. van Strien. One-dimensional
dynamics.  Springer, 1993.

\bibitem[NS]{NS} T. Nowicki, D. Sands.  Non-uniform hyperbolicity and
universal bounds for $S$-unimodal maps.  Invent. Math.  132  (1998),  no. 3,
633--680.

\bibitem[Pa]{Pa} J. Palis. A global view of dynamics and a Conjecture of the
denseness of finitude of attractors. 
Asterisque,
v. 261 (2000), 335 - 348. 


\bibitem[T1]{T} M. Tsujii.  Positive Lyapunov exponents in families of
one-dimensional maps.  Invent. Math. 111. 113-137, (1993).

\bibitem[T2]{T2} M. Tsujii.  Small random perturbations of one dimensional
dynamical systems and Margulis-Pesin entropy formula.
Random \& Comput. Dynamics. Vol.1 No.1 59-89, (1992).



\bibitem[Y]{Y} L.-S. Young.  Decay of correlations for certain
quadratic maps.  Comm. Math. Phys., 146 (1992), 123-138.

\end{thebibliography}
\end{document}